# Fast Power Series Solution of Large 3-D Electrodynamic Integral Equation for PEC Scatterers

Yoginder Kumar Negi[1], N. Balakrishnan[1], and Sadasiva M. Rao[2]

[1]Supercomputer Education Research Centre,
Indian Institute of Science, Bangalore, India
yknegi@gmail.com, balki@iisc.ac.in

[2]Naval Research Laboratory
Washington DC 20375, USA
sadasiva.rao@nrl.navy.mil

*Abstract*—This paper presents a new fast power series solution method to solve the Hierarchal Method of Moment (MoM) matrix for a large complex, perfectly electric conducting (PEC) 3D structures. The proposed power series solution converges in just two (2) iterations which is faster than the conventional fast solver–based iterative solution. The method is purely algebraic in nature and, as such applicable to existing conventional methods. The method uses regular fast solver Hierarchal Matrix (H-Matrix) and can also be applied to Multilevel Fast Multipole Method Algorithm (MLFMA). In the proposed method, we use the scaling of the symmetric near-field matrix to develop a diagonally dominant overall matrix to enable a power series solution. Left and right block scaling coefficients are required for scaling near-field blocks to diagonal blocks using Schur's complement method. However, only the right-hand scaling coefficients are computed for symmetric near-field matrix leading to saving of computation time and memory. Due to symmetric property, the left side-block scaling coefficients are just the transpose of the right-scaling blocks. Next, the near-field blocks are replaced by scaled near-field diagonal blocks. Now the scaled near-field blocks in combination with far-field and scaling coefficients are subjected to power series solution terminating after only two terms. As all the operations are performed on the near-field blocks, the complexity of scaling coefficient computation is retained as $O(N)$. The power series solution only involves the matrix-vector product of the far-field, scaling coefficients blocks, and inverse of scaled near-field blocks. Hence, the solution cost remains $O(NlogN)$. Several numerical results are presented to validate the efficiency and robustness of the proposed numerical method.

## I. INTRODUCTION

The integral equation-based Method of Moments (MoM) [1] is one of the popular methods for solving complex 3D electromagnetics problems. A few of the problems include scattering, radiation, EMI-EMC, etc. Compared to differential equation-based methods like Finite Difference Time Domain (FDTD) [2] and Finite Element Method (FEM) [3], the integral equation-based method results in fewer unknowns and is more stable and well-conditioned. With the recent advancement of computer speed and memory, the need for solving large size and complex problems in electromagnetics has increased rapidly. Conventional MoM is limited by quadratic memory storage, quadratic matrix fill time, and cubic solution time, which limits the application of MoM to resonance-size problems. To mitigate this issue and exploit the robustness of MoM, in the last few decades, more researchers have proposed fast solvers with $O(N \log N)$ matrix fill time and solution time, where $N$ represents the matrix size. A few of the popular fast solver methods for solving complex electrodynamic problems are Fast Fourier Transform (FFT), Multilevel Fast Multipole Algorithm (MLFMA) [4], IE-QR [5,6], Adaptive Cross Approximation (ACA) [7], etc. Most of these methods rely on the analytical or numerical matrix compression and fast matrix-vector product for solution leading to $O(NlogN)$ matrix fill and matrix-vector product cost. For multiple Right-Hand Side (RHS) problems like Mono-static Radar Cross Section (RCS) and Multi-port network with $N_{rhs}$ ports, and with $N_{itr}$ iterations for each solution, the total cost of the solution will be $O(N_{rhs}N_{itr}NlogN)$. Each solution of RHS is iteration dependent, and the iterations for desired tolerance depends on the condition number of the matrix. It is well known that an ill-conditioned matrix will lead to a high number of iterations, thus increasing the overall solution time. To improve the condition number of the matrix, researchers have suggested various types of matrix preconditioning methods like incomplete LU factorized ILUT [8], diagonal block-based Null-Field [9, 10] and Schur complement [11, 12]. But the effectiveness of these preconditioners is limited by precondition computation time and condition number improvement of the entire matrix. In contrast, the direct solver has an edge over the iterative solver, giving a solution in fixed single forward and backward solution operation for each RHS. However, the high-cost factorization, forward and back substitution limit the application of a direct solver for a large-size matrix. Recently, there is more inclination toward the development of fast direct solvers and various methods based on $H$ [13 -16].

In this work, we propose a method on par with a fast direct solver using the power series method, which converges in fixed 2 iterations. Recently, S. M. Rao and Michael S. Kluskens in [17] has proposed a method to solve the electromagnetic MoM matrix for electrically large conducting bodies by applying the power series method. The procedure involves the computation of the MoM matrix and dividing the matrix into relatively large subsections. The mutual coupling between a given subsection and related nearby subsections is transformed into self-coupling. The resulting current distribution is obtained by developing a power series solution. The power series method is suitable for solving bi-static and mono-static problems. The present work adopts the central idea of [17] and improves upon it to obtain a much faster solution.

In the present work, the procedure presented in [17] is modified using several essential steps. These steps include scaling near-field block matrices to diagonal block matrix and converting the scaled near-field in conjunction with far-field blocks to a power series format. Further, the diagonalization cost of computation and storage is reduced by using symmetric near-field blocks and adopting Adaptive Cross Approximation (ACA) for the far-field blocks [7]. As the diagonalization operation includes only near-field, the power series computation cost remains $O(N)$. The overall solution includes the matrix-vector product of the compressed far-field blocks and near-field blocks, retaining the overall solution cost as $O(NlogN)$. Further, extensive numerical experimentation shows that the proposed power series method converges in just two iterations. The present procedure is faster and efficiently applicable to large complex practical problems.

The paper is organized as follows: in section II, a brief description of multi-level CFIE H-Matrix 3D full-wave MoM is presented. Improved re-compressed ACA is used for matrix compression of the far-field blocks. In section III, the proposed power series format conversion from H-Matrix is presented along with the convergence criteria for the series. In section IV, complexity analysis of power series computation and memory cost is presented, and in section V, the efficiency and accuracy of the proposed power series solution are presented. Section VI concludes the paper.

## II. FAST H-MATRIX METHOD OF MOMENT

The 3D-electrodynamic problem for a PEC body can be solved using the Electric Field Integral Equation (EFIE), Magnetic Field Integral Equation (MFIE), or a Combined Field Integral equation (CFIE). The governing equation for EFIE states that the total electric field $E_{total}$ for a conducting 3D object is a combination of the incident field $E_{inc}$ and the scattered field, $E_{scatt}$

$$\boldsymbol{E}_{total} = \boldsymbol{E}_{inc} + \boldsymbol{E}_{scatt} \quad (1)$$

Applying the boundary condition for PEC surfaces, we have

$$\boldsymbol{E}_{inc} = j\omega\mu \int \int \mathbf{J}(\boldsymbol{r}') \, G(\boldsymbol{r},\boldsymbol{r}') \, ds'ds \\ + \frac{j}{\omega\epsilon} \int \int \boldsymbol{\rho}(\boldsymbol{r}')G(\boldsymbol{r},\boldsymbol{r}') \, ds'ds \quad (2)$$

where $\mathbf{J}(\boldsymbol{r}')$ and $\boldsymbol{\rho}(\boldsymbol{r}')$ represent the current density and charge density on the surface respectively, $\mu$ and $\epsilon$ represent the permeability and permittivity of the background material, $\omega$ is the angular frequency. In equation (2) $G$ is free-space Green's function and is given as

$$G(\boldsymbol{r},\boldsymbol{r}') = \frac{e^{jkR}}{4\pi R} \quad (3)$$

where $k$ is the wave-number, $\boldsymbol{r}$ and $\boldsymbol{r}'$ represent observer and source points and distance $\boldsymbol{R} = |\boldsymbol{r} - \boldsymbol{r}'|$. Integration is performed over observation surface $s'$ and source surface $s$.

Similarly, for MFIE [18], the boundary condition states that tangential total magnetic field $(\hat{\boldsymbol{n}} \times \boldsymbol{H}_{total})$ over a surface is equivalent to the electric current $(\boldsymbol{J}(\boldsymbol{r}))$ over the surface, as

$$\hat{\boldsymbol{n}} \times \boldsymbol{H}_{total} = \mathbf{J}(\boldsymbol{r}) \quad (4)$$

where, $\hat{\boldsymbol{n}}$ is a unit outward normal to the closed scattering surface. Now, the total magnetic field is a sum of incidence $(\boldsymbol{H}_i)$ and scattered magnetic field $(\boldsymbol{H}_{total})$

$$\hat{\boldsymbol{n}} \times (\boldsymbol{H}_{inc} + \boldsymbol{H}_{scatt}) = \mathbf{J}(\boldsymbol{r}) \quad (5)$$

Equation (5) can be extended further as

$$\hat{n} \times H_{inc} = \frac{J(r)}{2} - \hat{n} \times \int_{s'} J(r') \times \nabla' G(r,r') \, ds' \quad (6)$$

In the above equation $\nabla' G$ can be further simplified as

$$\nabla' G(r,r') = \left(jk + \frac{1}{R}\right) G(r,r') \hat{R} \quad (7)$$

Note that the EFIE is applicable for both open and closed surfaces, whereas MFIE is only applicable for a closed surface geometry. Combining EFIE and MFIE gives Combined Field Integral Equation (CFIE), given by

$$CFIE = \alpha EFIE + Z_o(1-\alpha) MFIE \quad (8)$$

where $\alpha$ is a control parameter to control the contribution of EFIE and MFIE, ranging from 0 to 1 and $Z_o$ is the free space impedance. The primary advantage of CFIE is that it is robust and generates a stable solution at internal resonances of the closed body and a well-conditioned MoM matrix. For an open surface, $\alpha$ is taken as 1, and for closed surfaces, $\alpha$ is taken as 0.5.

Current and charge density in surface integral equations EFIE and MFIE is modeled by RWG basis function [19], and Galerkin testing strategy is employed for MoM matrix computation. The final combined CFIE matrix is given as:

$$[Z] x = b \quad (9)$$

where $[Z]$ is a dense matrix of size $N \times N$ and $x$ and $b$ are unknown and known vectors of size $N \times 1$.

The CFIE dense matrix in equation (9) presents a time and memory bottleneck, with $O(N^2)$ memory, $O(N^2)$ matrix fill time, $O(N^3)$ for a direct solution and, $O(N_{itr}N^2)$ for iterative solver with $N_{itr}$ iterations. The iterative solution of the MoM matrix can be accelerated by exploiting the compressibility of the far-field sub-matrices. Compression also expedites the cost of matrix fill time and matrix-vector multiplication time in an iterative solver. Compression can be done analytically, like in the case of MLFMA or algebraically using IE-QR or ACA. Due to the kernel-independent property of the algebraic compression method, recently, ACA has gained popularity among researchers for the development of fast solvers. These compression methods can be applied in conjunction with binary-tree–based multi-level Hierarchal Matrix (H-Matrix). In the binary-tree decomposed 3D geometry, the matrix compression is applied for block interaction at each level, satisfying the admissibility condition given below

$$\eta \, dis(\Omega_t, \Omega_s) \geq \min(dia(\Omega_t), dia(\Omega_s)) \quad (10)$$

The admissibility condition states that the minimum of the block diameter of the test block ($\Omega_t$) and source block ($\Omega_s$) should be less than or equal to the admissibility constant ($\eta$) times the distance between the test and source blocks. The binary-tree partition of the geometry is carried out until the block size is greater than or equal to $0.5\lambda$. The criteria for binary-tree truncation are discussed in section V. At the leaf level, the block interaction not satisfying the admissibility condition is considered as a near-field interaction.

In this work, the re-compressed ACA method [20, 21] is employed for the computation of multi-level H-Matrix. For the $m \times n$ rectangular sub-matrix $Z_{sub}^{m \times n}$ representing the coupling between two well-separated groups of $m$ observer bases and $n$ source bases, the ACA algorithm aims to approximate $Z_{sub}^{m \times n}$ by $A^{m \times r}$ and $B^{r \times n}$ such that:

$$Z_{sub}^{m \times n} \approx A^{m \times r} \times B^{r \times n} \quad (11)$$

where $r$ is the effective rank of the matrix $Z_{sub}^{m \times n}$ such that $r \ll \min(m,n)$, $A^{m \times r}$ and $B^{r \times n}$ are two low rank dense rectangular matrices, satisfying the accuracy condition

$$\|Z_{sub}^{m \times n} - A^{m \times r} \times B^{r \times n}\| \leq \varepsilon \|Z_{sub}^{m \times n}\| \quad (12)$$

For a given tolerance $\varepsilon$ $\|.\|$ refers to the matrix Frobenius norm. Traditional ACA-based methods suffer from higher rank and error for the desired tolerance. To mitigate this, a re-compression scheme is suggested in [20, 21]. The compression cost for each sub-matrix is given by $r^2(m+n)$ and the storage and matrix-vector product cost by $r(m+n)$. The multi-level binary-tree matrix decomposed H-Matrix method leads to $O(NlogN)$ matrix fill and matrix-vector product time for each iteration. However, the final solution cost is highly iteration dependent as for $N_{itr}$ iterations, the solution cost scales to $O(N_{itr}NlogN)$. Further, for the case of multiple RHS with $N_{rhs}$ vectors, the solution cost scales to $O(N_{rhs}N_{itr}NlogN)$. To mitigate this iteration-dependent solution, we have proposed a power series–based iterative solver which converges in fixed 2 iterations maintaining the optimum cost of a fast solver.

### III. POWER SERIES SOLUTION

In this section, we present a new fast power series solution method for solving large H-Matrix. The method follows Schur's complement procedure used for matrix diagonalization and precondition computation [11, 12]. Power series is an infinite series and a well-known method for solving an ordinary and partial differential equation. The advantage of the power series is that it converges in a very small region. Since a conventional MoM matrix cannot be applied for the power series solution, we propose a method to convert to the H-matrix format. We note that the H-matrix is suitable for the power series solution, as discussed further in the section. H-Matrix is a combination of the far-field and near-field matrices. In the following subsections, the details of the conversion process are presented.

## A. Preparing for Power Series Computation

As a first step, the geometry is divided into blocks based on the same binary tree as used in the compression algorithm. In a multi-level H-Matrix compression scheme, as described in section II, the MoM matrix $[Z]$ in equation (9) can be represented as a combination of the near-field $[Z_N]$ matrix at the leaf level and the compressed far-field $[Z_F]$ matrix obtained at multiple levels of far-interaction.

$$[Z]x = [Z_N + Z_F]x = b \quad (13)$$

where $x$ represents the unknown coefficient vector, $b$ is the excitation vector. To maintain the optimum cost, symmetric near-field matrix is used. To explain the procedure, a leaf-level cube structure comprised of four cubes, as shown in Figure 2, is considered.

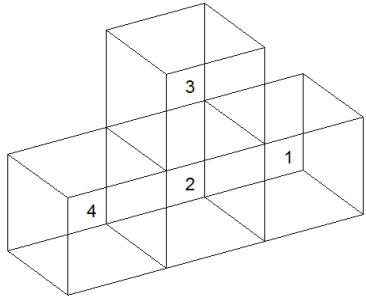

Figure 1 Representative leaf-level cubes for illustration of Schur's process.

For Figure 1, block interaction between 1 and 4 forms the far-field interaction and rest form the near-field interaction. So near-field block matrix $[Z_N]$ in the case of Figure 1 is given as:

$$[Z_N] = \begin{bmatrix} Z_{11} & Z_{12} & Z_{13} & 0 \\ Z_{21} & Z_{22} & Z_{23} & Z_{24} \\ Z_{31} & Z_{32} & Z_{33} & Z_{34} \\ 0 & Z_{42} & Z_{43} & Z_{44} \end{bmatrix} \quad (14)$$

Now, the near-field can be scaled completely to a diagonal block by using the left- and right-hand scaling coefficients. Right scaling coefficient $[\alpha_1]$ for scaling the first-row blocks $Z_{11}$ and $Z_{13}$ can be represented as:

$$[\alpha_1] = \begin{bmatrix} I_{11} & \alpha_{12} & \alpha_{13} & 0 \\ 0 & I_{22} & 0 & 0 \\ 0 & 0 & I_{33} & 0 \\ 0 & 0 & 0 & I_{44} \end{bmatrix} \quad (15)$$

where, $I_{11}$, $I_{22}$, $I_{33}$, and $I_{44}$ are the identity block matrices. In the null-field method [9, 10], either left- or right-hand scaling is performed for the near-field scaling to diagonal blocks ignoring the fill-in blocks. In contrast, in the proposed Schur's complement method, both left- and right-hand scaling are performed simultaneously, considering fill-in blocks for complete far-field scaling to diagonal blocks. The fill-ins in the scaling matrix are described in the following sub-section for row scaling, the values in $[\alpha_1]$ can be given as:

$$[\alpha_1] = \begin{bmatrix} I_{11} & -Z_{11}^{-1}Z_{12} & -Z_{11}^{-1}Z_{13} & 0 \\ 0 & I_{22} & 0 & 0 \\ 0 & 0 & I_{33} & 0 \\ 0 & 0 & 0 & I_{44} \end{bmatrix} \quad (16)$$

Similarly, for the complete scaling of column blocks $Z_{21}$ and $Z_{31}$ the left scaling coefficient $[\alpha'_1]$ is used and $[\alpha'_1]$ can be given as:

$$[\alpha'_1] = \begin{bmatrix} I_{11} & 0 & 0 & 0 \\ \alpha'_{12} & I_{22} & 0 & 0 \\ \alpha'_{13} & 0 & I_{33} & 0 \\ 0 & 0 & 0 & I_{44} \end{bmatrix} \quad (17)$$

$$[\alpha'_1] = \begin{bmatrix} I_{11} & 0 & 0 & 0 \\ -Z_{21}Z_{11}^{-1} & I_{22} & 0 & 0 \\ -Z_{31}Z_{11}^{-1} & 0 & I_{33} & 0 \\ 0 & 0 & 0 & I_{44} \end{bmatrix} \quad (18)$$

Now, equations (16), (18), and (14) can be combined to scale the first row and column block of $[Z_N]$ to diagonal block and the system of the equation can be given as:

$$[\tilde{Z}_N^1] = \begin{bmatrix} I_{11} & 0 & 0 & 0 \\ \alpha'_{12} & I_{22} & 0 & 0 \\ \alpha'_{13} & 0 & I_{33} & 0 \\ 0 & 0 & 0 & I_{44} \end{bmatrix} \begin{bmatrix} Z_{11} & Z_{12} & Z_{13} & 0 \\ Z_{21} & Z_{22} & Z_{23} & Z_{24} \\ Z_{31} & Z_{32} & Z_{33} & Z_{34} \\ 0 & Z_{42} & Z_{43} & Z_{44} \end{bmatrix}$$
$$\times \begin{bmatrix} I_{11} & \alpha_{12} & \alpha_{13} & 0 \\ 0 & I_{22} & 0 & 0 \\ 0 & 0 & I_{33} & 0 \\ 0 & 0 & 0 & I_{44} \end{bmatrix} \quad (19)$$

Equation (19) can be represented as:

$$[\tilde{Z}_N^1] = [\alpha'_1][Z_N][\alpha_1] \quad (20)$$

Performing the block multiplication in equation (19), $[\tilde{Z}_N^1]$ can be represented as a block matrix form as:

$$[\tilde{Z}_N^1] =$$
$$\begin{bmatrix} Z_{11} & 0 & 0 & 0 \\ 0 & Z_{22}-Z_{21}Z_{11}^{-1}Z_{12} & Z_{23}-Z_{21}Z_{11}^{-1}Z_{13} & Z_{24} \\ 0 & Z_{32}-Z_{31}Z_{11}^{-1}Z_{12} & Z_{33}-Z_{31}Z_{11}^{-1}Z_{13} & Z_{34} \\ 0 & Z_{42} & Z_{43} & Z_{44} \end{bmatrix} \quad (21)$$

Equation (21) gives Schur's complement of the first block near-field matrix. Likewise, each row and column block can be scaled to form a diagonal block matrix and is of the form:

$$[\widetilde{Z}_N] = [\,\alpha'_3\,][\,\alpha'_2\,][\,\alpha'_1\,][Z_N][\alpha_1][\alpha_2]][\alpha_3] \quad (22)$$

$$[\widetilde{Z}_N] = \begin{bmatrix} Z_{11} & 0 & 0 & 0 \\ 0 & \widetilde{Z}_{22} & 0 & 0 \\ 0 & 0 & \widetilde{Z}_{33} & 0 \\ 0 & 0 & 0 & \widetilde{Z}_{44} \end{bmatrix} \quad (23)$$

Equation (23) gives the complete diagonal form of the near-field matrix. For solving the complete system of equations with left- and right-hand scaling coefficients, the final system of equation (9) can be represented as:

$$[\,\alpha'_3\,][\,\alpha'_2\,][\,\alpha'_1\,][Z][\alpha_1][\alpha_2]][\alpha_3][\widetilde{x}] = [\widetilde{b}] \quad (24)$$

Now, $[b]$ and $[x]$ in equation (9) can be extracted by

$$[\widetilde{b}] = [\,\alpha_3^T\,][\,\alpha_2^T\,][\,\alpha_1^T\,][b] \quad (25)$$

$$[x] = [\alpha_1][\alpha_2]][\alpha_3][\widetilde{x}] \quad (26)$$

Equation (24) can be defined as the sum of the near and far-field as in equation (13) and is given as:

$$[\,\alpha'_3\,][\,\alpha'_2\,][\,\alpha'_1\,][Z_N + Z_F][\alpha_1][\alpha_2][\alpha_3][\widetilde{x}] = [\widetilde{b}] \quad (27)$$

where, $[Z_F]$ is the far-field compressed ACA matrix blocks and $[Z_N]$ is the dense near-field block matrices. Equation (24) can be further simplified by as:

$$[\alpha'_3][\,\alpha'_2\,][\alpha'_1][Z_N][\alpha_1][\alpha_2][\alpha_3][\widetilde{x}] + \\ [\alpha'_3][\alpha'_2][\alpha'_1][Z_F][\alpha_1][\alpha_2][\alpha_3][\widetilde{x}] = [\widetilde{b}] \quad (28)$$

The first part of the above equation represents the block diagonal near-field as in equation (23), and then the equation can be further simplified as:

$$[\widetilde{Z}_N][\widetilde{x}] + [\alpha'_3][\alpha'_2][\alpha'_1][Z_F][\alpha_1][\alpha_2][\alpha_3][\widetilde{x}] = [\widetilde{b}] \quad (29)$$

where, $[\widetilde{Z}_N]$ is a scaled near-field block diagonal matrix. Due to the symmetric property of the near-field matrix, we only need to compute and store right-hand scaling coefficients $[\alpha_1]$, $[\alpha_2]$, and $[\alpha_3]$ since left-hand scaling coefficients $[\alpha'_1]$, $[\alpha'_2]$, and $[\alpha'_3]$ are mere transpose of right-hand coefficients. Figure 2. below shows the sparsity pattern of the near-field and right-hand scaling coefficient matrix for $5\lambda \times 5\lambda$ plate. Sloan's graph ordering is used to reduce the fill-in and computation cost for the scaling coefficient, as suggested in [11].

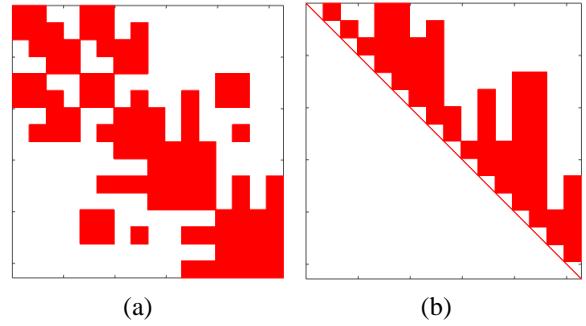

(a)          (b)
Figure 2. Sparsity pattern of the (a) near-field matrix, (b) right-hand scaling coefficient.

Now equation (29) can be converted to power series solution format by moving the scaled diagonal near-field to the right-hand side of the equation leads to

$$\left[[I] + [\widetilde{Z}_N]^{-1}[\alpha'_3][\alpha'_2][\alpha'_1][Z_F][\alpha_1][\alpha_2][\alpha_3]\right][\widetilde{x}] = \\ [\widetilde{Z}_N]^{-1}[\widetilde{b}] \quad (30)$$

Each term in equation (30) can be represented as

$$[U] = [\widetilde{Z}_N]^{-1}[\alpha'_3][\alpha'_2][\alpha'_1][Z_F][\alpha_1][\alpha_2][\alpha_3] \quad (31)$$

$$[b_o] = [\widetilde{Z}_N]^{-1}[\widetilde{b}] \quad (32)$$

Substituting equation (31) and (32) in equation (30) leads to

$$[[I] + [U]][\widetilde{x}] = [b_o] \quad (33)$$

$$[\widetilde{x}] = [[I] + [U]]^{-1}[b_o] \quad (34)$$

Equation (34) can be solved using power series solution

$$[\widetilde{x}] = [[I] - [U] + [U]^2 - [U]^3 + \ldots\ldots][b_o] \quad (35)$$

$$[\widetilde{x}] = [b_o] - [Ub_o] + [U[Ub_o]] - [U[U[Ub_o]]] + \ldots \quad (36)$$

Equation (36) shows that solving equation (34) is an iterative matrix-vector product of the inverse of the scaled diagonal block near-field, scaling coefficients and far-field blocks.

**B. Convergence**

Power series always converges in the radius of convergence. The necessary and sufficient condition for the power series equation (35) to converge in the radius of convergence is the Frobenius norm, $\|U\| \leq 1$ in

equation (31). Defining

$$[\tilde{Z}_F] = [\alpha'_3][\alpha'_2][\alpha'_1][Z_F][\alpha_1][\alpha_2]][\alpha_3] \quad (37)$$

Equation (28) can be re-written as

$$\left[[I] + [\tilde{Z}_N]^{-1}[\tilde{Z}_F]\right][\tilde{x}] = [\tilde{Z}_N]^{-1}[\tilde{b}] \quad (38)$$

To achieve the convergence, we can enforce the condition $\left\|[\tilde{Z}_N]^{-1}\right\| \cdot \left\|[\tilde{Z}_F]\right\| \leq 1$ for the final solution of the power series. Since the process of norm computation for a large matrix is compute-intensive, alternatively, one can adopt the following procedure. We note that

$$\left\|[\tilde{Z}_N]^{-1}\right\| = \frac{\left\|[\tilde{Z}_N]^{-1}\right\| \cdot \left\|[\tilde{Z}_N][b_o]\right\|}{\left\|[\tilde{Z}_N][b_o]\right\|}$$

$$\leq \frac{\left\|[\tilde{Z}_N]^{-1}\right\| \cdot \left\|[\tilde{Z}_N]\right\| \cdot \left\|[b_o]\right\|}{\left\|[\tilde{b}]\right\|}$$

$$= k_{nf} \frac{\left\|[b_o]\right\|}{\left\|[\tilde{b}]\right\|} \quad (39)$$

where, $k_{nf} = \left\|[\tilde{Z}_N]^{-1}\right\| \cdot \left\|[\tilde{Z}_N]\right\|$ represents the condition number of $[\tilde{Z}_N]$.

Next, we define $[b_e] = [\tilde{Z}_F][b_o]$ and we have

$$\left\|[\tilde{Z}_F]\right\| = \frac{\left\|[\tilde{Z}_F]\right\| \cdot \left\|[[\tilde{Z}_F]^{-1}][b_e]\right\|}{\left\|[\tilde{Z}_F]^{-1}[b_e]\right\|}$$

$$\leq \frac{\left\|[\tilde{Z}_F]\right\| \cdot \left\|[[\tilde{Z}_F]^{-1}]\right\| \cdot \left\|[b_e]\right\|}{\left\|[b_o]\right\|}$$

$$= k_{ff} \frac{\left\|[b_e]\right\|}{\left\|[b_o]\right\|} \quad (40)$$

where, $k_{ff} = \left\|[\tilde{Z}_F]^{-1}\right\| \cdot \left\|[\tilde{Z}_F]\right\|$ represents the condition number of $[\tilde{Z}_F]$.

Combining equations (39) and (40), we have

$$\left\|[\tilde{Z}_N]^{-1}\right\| \cdot \left\|[\tilde{Z}_F]\right\| \leq k_{nf} k_{ff} \frac{\left\|[b_e]\right\|}{\left\|[\tilde{b}]\right\|} \quad (41)$$

To satisfy the condition $\|U\| \leq 1$ we must ensure that

$$\frac{\left\|[b_e]\right\|}{\left\|[\tilde{b}]\right\|} \leq \frac{1}{k_{nf} k_{ff}} \quad (42)$$

$$\left\|[\tilde{Z}_N]^{-1}[\tilde{Z}_F]\right\| \leq \frac{1}{k_{nf} k_{ff}} \quad (43)$$

Now, let us consider each right-hand side term in equation (36) is represented as the sum of iteration terms $it_o, it_1, it_2 \ldots it_n$ leading to

$$[\tilde{x}] = it_o - it_1 + it_2 \ldots (-1)^n it_n \quad (44)$$

For the convergence test, we can check the iteration norm ratio as:

$$\frac{\|it_n\|}{\|it_{n-1}\|} = \|U\| = \left\|[\tilde{Z}_N]^{-1}[\tilde{Z}_F]\right\| \leq \frac{1}{k_{nf} k_{ff}} \quad (45)$$

Note that it is not really necessary to compute the condition numbers $k_{nf}$ and $k_{ff}$ but ensure that the fraction $\frac{\|it_n\|}{\|it_{n-1}\|}$ is a small number. Our various numerical experiments suggested that this number must be less than $1e^{-1}$. It is because the described numerical implementation ensures that matrices $[\tilde{Z}_N]$ and $[\tilde{Z}_F]$ are well-conditioned matrices. Obviously, if the $\frac{\|it_n\|}{\|it_{n-1}\|}$ is not less than the empirical value, then the solution may diverge.

## IV. COMPLEXITY ANALYSIS

In this section, the linear order complexity for power series set-up time is presented. For the complexity analysis, a uniform distribution of $N$ RWG bases in 3D grouped in a cube, and following a multi-level binary-tree decomposition, each cube is recursively subdivided into two cubes starting from level 0 to level $L$. Therefore, at the lowest level, there are $2^L$ leaf-level cubes. Assuming a uniform distribution, the number of basis functions in each leaf-level cube is $\frac{N}{2^L}$. Also, following the theory of most fast solver algorithms, it can be shown that for optimal efficiency of matrix storage and matrix-vector product cost $L = \log_2 N$.

### A. Computation Cost

The power series set-up cost includes scaling the near-field matrix to diagonal format and arranged scaled near-field, far-field, and near-field scaling coefficients to a power series format. Near-field scaling cost is the high cost of power series set-up. The near-field scaling consists of computation right scaling coefficients $[\alpha]$ as in equation (16). Due to the symmetric property, left scaling coefficients are just the transpose of right-hand scaling coefficients. Right-hand scaling coefficient computation cost can be represented as $C_1$. The second cost includes the scaling of the near-field to the diagonal block form by $[\alpha'][Z_N][\alpha]$ operation. For each row and column block in equations (18) and (20), this cost can be represented as $C_2$. Therefore, the total cost can be summed up as:

$$C_{TOTAL} = C_1 + C_2 \qquad (46)$$

1. **Scaling coefficient computation cost**

For the scaling coefficient computation, the high cost includes the inversion ($C_{MI}$) cost for diagonal block and the solving the inverse ($C_{SOL}$) for the row and column block far-fields as in equations (14) and (16). Therefore $C_{SCC}$ can be further be divided as the summation of inversion and solution cost as:

$$C_1 = C_{MI} + C_{SOL} \qquad (47)$$

Inversion cost includes the single matrix inversion of a diagonal block for scaling near-field of each row and column block. Therefore, the matrix inversion cost of one matrix of a diagonal block at leaf level is given as:

$$C_{MI}^1 = k_1 \times \left[\frac{N}{2^L}\right]^3 \qquad (48)$$

where, $k_1$ is a constant, the total cost for matrix inversion for the leaf level blocks is given by

$$C_{MI} = \sum_{i=1}^{2^L} C_{MI}^i = k_1 \times \left[\frac{N}{2^L}\right]^3 \times 2^L \qquad (49)$$

$$C_{MI} = k_1 \times N = O(N) \qquad (50)$$

For the computation of the scaling coefficient, the inverted matrix has to be solved for all row and column near-field blocks. Cost of matrix solution for one block at leaf level can be given as:

$$C_{SOL}^1 = k_2 \times \left[\frac{N}{2^L}\right]^2 \times \left[\frac{N}{2^L}\right] \qquad (51)$$

where, $k_2$ is a constant. For a 3D structure, each block is surrounded by 26 near-field blocks. Therefore, the matrix solution cost for each row is given by:

$$C_{SOL}^{1R} = k_2 \times \left[\frac{N}{2^L}\right]^2 \times 26 \times \left[\frac{N}{2^L}\right] \qquad (52)$$

The total cost of the matrix solution at the leaf level blocks is the summation of the cost of each row and is given as:

$$C_{SOL} = \sum_{i=1}^{2^L} C_{SOL}^{iR} = k_2 \times \left[\frac{N}{2^L}\right]^2 \times 26 \times \left[\frac{N}{2^L}\right] \times 2^L \qquad (53)$$

$$C_{SOL} = k_2 \times 26 \times N = O(N) \qquad (54)$$

2. **Near-Field scaling cost**

For converting near-field to a diagonal block matrix format, the near-field block matrix has to be multiplied by the left and right scaling coefficients as given in equations (19) and (20). From equation (21), the multiplication involves block near-field matrix and column block matrix

$$c_S^1 = k_3 \times \left[26\frac{N}{2^L} \times 26\frac{N}{2^L}\right] \times \left[\frac{N}{2^L} \times 26\frac{N}{2^L}\right] \qquad (55)$$

where $k_3$ is a constant. Therefore, the total cost of scaling the near-field blocks with the right-hand scaling coefficients blocks is

$$C_S = \sum_{i=1}^{2^L} C_S^i \approx O(N) \qquad (56)$$

Equations (50), (54), and (56) show that the cost of the power series computation, for a uniform 3D distributed basis function, is $O(N)$. As the storage of scaling coefficient and scaled near-field block is half the near-field matrix size, the memory cost scales to $O(N)$. The $O(N)$ complexity of power series computation and memory is experimentally shown in Figure 3 for the increasing number of unknowns and the size of a sphere.

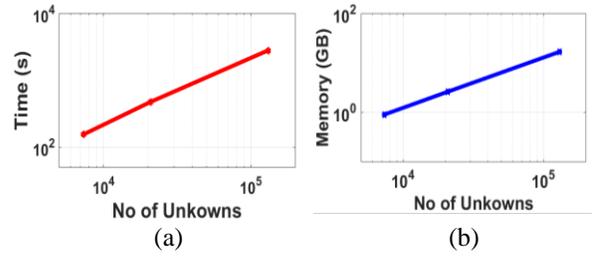

Figure 3 (a) Power series set-up time, (b) Memory in GB for scaling coefficient and scaled near-field with increasing unknowns.

## V. NUMERICAL RESULTS

In this section, we show the binary-tree truncation criteria and solution complexity for the power series solution. The accuracy and efficiency are shown by RCS comparison for different geometries. All the computations were carried out for double-precision data type on the system with 128 GB memory and Intel Xeon E5-2670 processor. The comparisons are made for an open and closed structure.

### A. Binary-Tree Truncation

For binary-tree truncation, we tested the accuracy of the power series solution vector from $2\lambda$ sphere of 20,802 unknowns for varying binary-tree leaf-level sizes. For accuracy check, Frobenius norm error was calculated between solution vectors from direct LU factorized solution $SOL_{dir}$ and power series solution $SOL_{ps}$ by computing $\|SOL_{dir} - SOL_{ps}\|/\|SOL_{dir}\|$. The error plot is shown in Figure 4 below

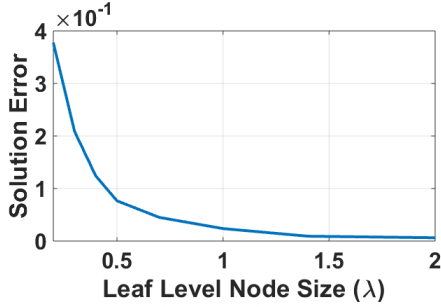

Figure 4 Power Series Solution error with increasing leaf level node size.

It can be observed from Figure 4 above that for the leaf level truncation size greater than $0.5\,\lambda$ we get the desired accuracy. Hence, for all simulations shown in this work, the binary-tree leaf node is truncated for a size greater than 0.5.

### B. Solution Complexity

In Figure 5, we demonstrate that the proposed power series method retains the $O(NlogN)$ solution complexity of H-Matrix. The experiment is carried out for sphere meshed with increasing sphere size and unknowns.

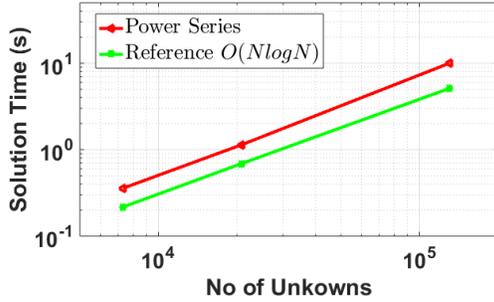

Figure 5 Power Series solution time with increasing unknown.

### C. Accuracy and Efficiency

In this subsection, to validate the accuracy of the proposed method, we have compared the RCS results of different geometries with the analytical results. Also, to demonstrate the efficiency of the proposed method, solution time is compared with regular iterative and preconditioned [11] iterative solutions. For all comparative case studies, iterative solver GMRES with an error tolerance of 1e-6 is considered.

#### 1. Bi-Static RCS of a sphere

As a first example, we consider bi-static RCS of a $5\lambda - radius$ sphere discretized with a $\lambda/10$ mesh resulting in 130,293 unknowns. The solution from the method described in this work is compared with the Mie series analytical solution. Figure 6 shows the agreement of bi-static RCS from the present method with the Mie series. RCS is computed for observation angle $\theta = 0^o$ to $180^o$ for $\phi = 0^o$ with VV polarized plane wave incident at $\theta = 0^o$ and $\phi = 0^o$. We note excellent agreement between the two results.

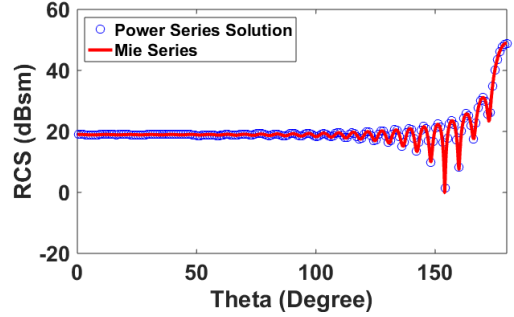

Figure 6 Bi-static RCS of $5\lambda$ sphere for observation angles $\theta = 0^o$ to $180^o$, $\phi = 0^o$ and VV polarized plane wave incident at $\theta = 0^o, \phi = 0^o$.

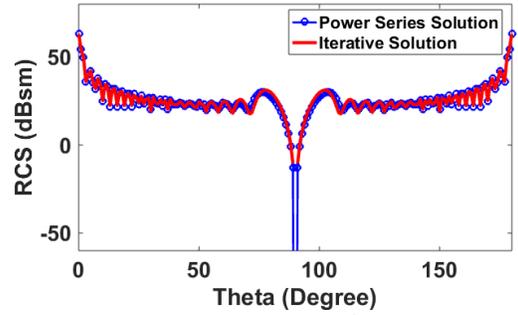

Figure 7 Mono-static RCS of $20\lambda$ square plate with VV polarized incident plane wave and observation angles at $\theta = 0^o$ to $180^o$, $\phi = 0^o$

#### 2. Mono-Static RCS of a Square Plate

To show the accuracy and efficiency of the proposed power series method, an open structure, a square plate of $20\lambda$ size, is considered. The solution is obtained using EFIE only with $\alpha = 1$ in equation (8). Figure 7 shows mono-static RCS of the plate using the present method and compared with a conventional iterative method. The plate is discretized with $\lambda/10$ element size for 119,600 unknowns. The square plate is located in the XY plane and illuminated by a plane wave incident with $\theta$ varying from $0\ to\ 180^o$ and $\phi = 0^o$. It can be observed from the figure that there is a very favorable agreement of RCS between the two methods.

TABLE I
SOLUTION TIME OF SQUARE PLATE

| Method | Setup Time (H) | Solution Time (H) |
|---|---|---|
| Power Series Solver | 1.61 | 0.57 |
| Preconditioned Iterative Solver | 1.61 | 3.760 |
| Iterative Solver | -------- | 89.97 |

In Table 1, we present the comparison of solution time for the results presented in Figure 7. We show solution times for iterative solvers with and without preconditioning. It is evident that the present work is much more efficient than the other two cases.

### 3. Mono-Static RCS of a Cube

Next, we consider a conducting cube of 1 m meshed with $\lambda/10$ element size giving 45,975 unknowns. Since the scattering body is a closed structure, we use CFIE with $\alpha = 0.5$. The operating frequency is 1.3 GHz and compared with the iterative solution [20], as shown in Figure 8. Further, in Table II, we present a comparison of solution time for this example.

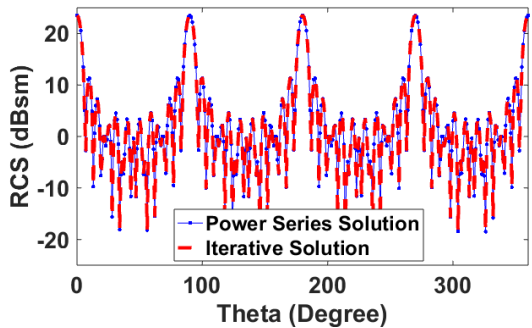

Figure 8 Mono-static RCS of 1 $m$ cube at 1.3 GHz for HH polarized plane wave incident and observation angles at $\theta = 0^o$ to $360^o$, $\phi = 0^o$.

TABLE II
SOLUTION TIME OF A CUBE

| Method | Setup Time (H) | Solution Time (H) |
|---|---|---|
| Power Series Solver | 1.32 | 0.551 |
| Preconditioned Iterative solver | 1.32 | 1.377 |
| Iterative solver | -------- | 7.640 |

It can be observed that the RCS result from the power series completely matches with the regular H-matrix iterative solver [22].

### 4. Mono-Static RCS of a Model Fighter Aircraft

As the last example, we consider the geometry of model fighter aircraft with length 14 m and wingspan 8 m. With $\lambda/10$ discretization of the geometry, the meshing scheme generates 93,819 unknowns. Figure 9 shows the computed mono-static RCS with $\alpha = 0.5$ in CFIE equation (8) at 300 MHz in the X-Y plane with VV polarized plane wave incident and observation angle at the nose to tail $\phi = 0^o$ to $180^o$ and $\theta = 90^o$. Mono-static RCS (180 RHS) solution time comparisons of power series and iterative solver with and without preconditioned are shown in Table III.

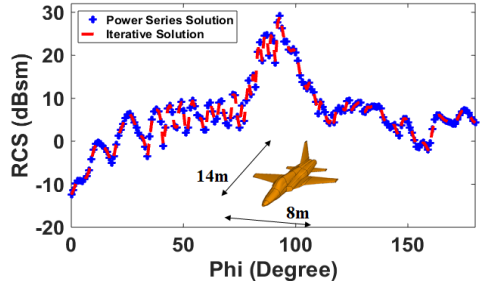

Figure 9 Mono-static RCS of model fighter aircraft at 300 MHz for VV polarized plane wave incident and observation angles at $\theta = 90^o$, $\phi = 0^o$ to $180^o$.

TABLE III
SOLUTION TIME OF A MODEL FIGHTER AIRCRAFT

| Method | Setup Time (H) | Solution Time (H) |
|---|---|---|
| Power Series Solver | 1.66 | 0.46 |
| Preconditioned Iterative solver | 1.66 | 1.95 |
| Iterative solver | -------- | 191.73 |

It can be observed from Figure 9 that for this complex geometry, RCS from the power series solution entirely agrees with the H-Matrix iterative solution and results in a much higher efficient solution.

## VI. CONCLUSION

In this work, we propose a new power series solution method for solving 3D MoM-based integral equations. It can be observed from the numerical experimentation that the proposed method is as accurate as of the conventional iterative H-Matrix solution. Also, the proposed power series method results in significantly lower solution time compared to regular iterative and preconditioned iterative solutions. The method is based on the near-field matrix operation, thus maintaining $O(N)$ complexity for computation and $O(N\log N)$ solution time. The solution converges in fixed 2 iterations. The proposed method is a kernel-independent algebraic method and can be applied to other acceleration algorithms like MLFMA.

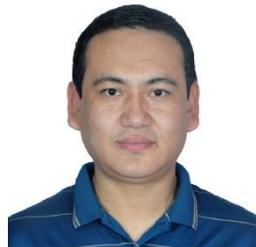

**Yoginder Kumar Negi** obtained the B.Tech degree in Electronics and Communication Engineering from Guru Gobind Singh Indraprastha University, New Delhi, India, in 2005, M.Tech degree in Microwave Electronics from Delhi University, New Delhi, India, in 2007 and the PhD degree in engineering from Indian Institute of Science (IISc), Bangalore, India, in 2018.

Dr Negi joined Supercomputer Education Research Center (SERC), IISc Bangalore in 2008 as a Scientific Officer. He is currently working as a Senior Scientific Officer in SERC IISc Bangalore. His current research interests include numerical electromagnetics, fast techniques for electromagnetic application, bio-electromagnetics, high-performance computing, and antenna design and analysis.
.

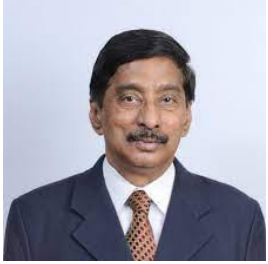

**B Narayanaswamy** received the B.E. degree (Hons.) in Electronics and Communication from the University of Madras, Chennai, India, in 1972, and the Ph.D. degree from the Indian Institute of Science, Bengaluru, India, in 1979.

He joined the Department of Aerospace Engineering, Indian Institute of Science, as an Assistant Professor, in 1981, where he became a Full Professor in 1991, served as the Associate Director, from 2005 to 2014, and is currently an INSA Senior Scientist at the Supercomputer Education and Research Centre. He has authored over 200 publications in the international journals and international conferences. His current research interests include numerical electromagnetics, high-performance computing and networks, polarimetric radars and aerospace electronic systems, information security, and digital library.

Dr. Narayanaswamy is a fellow of the World Academy of Sciences (TWAS), the National Academy of Science, the Indian Academy of Sciences, the Indian National Academy of Engineering, the National Academy of Sciences, and the Institution of Electronics and Telecommunication Engineers.

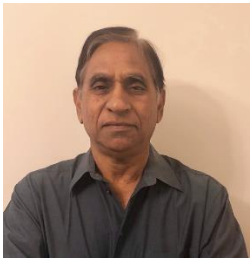

**Sadasiva M. Rao** obtained his Bachelors, Masters, and Doctoral degrees in electrical engineering from Osmania University, Hyderabad, India, Indian Institute of Science, Bangalore, India, and University of Mississippi, USA, in 1974, 1976, and 1980, respectively. He is well known in the electromagnetic engineering community and included in the Thomson Scientifics' *Highly Cited Researchers List*.

Dr. Rao has been teaching electromagnetic theory, communication systems, electrical circuits, and other related courses at the undergraduate and graduate level for the past 30 years at various institutions. At present, he is working at Naval Research Laboratories, USA. He published/presented over 200 papers in various journals/conferences. He is an elected Fellow of IEEE.